\documentclass[final,onefignum,onethmnum,oneeqnum]{siamltex}
\usepackage{amsmath}  
\usepackage{epsfig}

\newcommand{\dx}{{{\mit \Delta} x}}
\newcommand{\dt}{{{\mit \Delta} t}}

\newcommand{\e}{{\rm e}}

\newcommand{\nn}{\nonumber}
\newcommand{\la}{\label} 

\newcommand{\be}{\begin{equation}}
\newcommand{\ee}{\end{equation}}
\newcommand{\ba}{\begin{eqnarray}}
\newcommand{\ea}{\end{eqnarray}}
\newcommand{\bt}{\begin{theorem}}
\newcommand{\et}{\end{theorem}}
\newcommand{\bp}{\begin{proof}}
\newcommand{\ep}{\end{proof}}
\title{A theory of explicit finite-difference schemes}

\author{Siu A. Chin\thanks{Department of Physics and Astronomy, Texas A\&M University,
College Station, TX 77843, U.S.A. ({\tt chin@physics.tamu.edu}). }}

\begin{document}

\maketitle

\begin{abstract}
Conventional finite-difference schemes for solving partial differential
equations are based on approximating derivatives by finite-differences. 
In this work, an alternative theory is proposed which view finite-difference 
schemes as systematic ways of matching up to the operator solution of the 
partial differential equation. By completely abandon the idea of approximating 
derivatives directly, the theory provides 
a unified description of explicit finite-difference schemes for solving a 
general linear partial differential equation with constant coefficients to any 
time-marching order. As a result, the stability of the first-order 
algorithm for an entire class of linear equations can be determined all at once.
Because the method is based on solution-matching, it can also be used to derive
any order schemes for solving the general nonlinear advection equation.   
\end{abstract}

\begin{keywords} 
Finite-difference schemes, higher order methods, Burgers' equation.
\end{keywords}

\begin{AMS}
65M05, 65M06, 65M10
\end{AMS}

\pagestyle{myheadings}
\thispagestyle{plain}
\markboth{SIU A. CHIN}{A THEORY OF FINITE-DIFFERENCE SCHEMES}

\section {Introduction}
The most fundamental aspect of devising numerical algorithms for solving
partial differential equations is to derive finite-difference schemes for 
solving a general linear equation of the form
\be
\frac{\partial u}{\partial t}=\sum_{m=1}^M\!a_m\partial_x^mu,
\la{geq}
\ee																	  
with constant coefficients $a_m$. Conventionally, numerical schemes are obtained by
approximating the temporal and spatial derivatives of the equation by finite-differences. 
Such a direct use of finite-difference approximations produces a large collection
of seemingly unrelated and disparate finite-difference schemes which must
be analyzed one by one for stability and efficiency. There does not appear to be a
unifying theme that connects all such schemes. 
Moreover, if only explicit schemes are desired,
then discretizing the equation can only produce first-order algorithms,  
since the required grid values at multiple time steps can only be obtained implicitly.

To go beyond first-order, instead of approximating the
equation, one can approximate the formal operator solution to the
equation. In the case of (\ref{geq}), the solution is 
\be
u(x,\dt)=\exp\Bigl(\dt\sum_{m=1}^Ma_m\partial_x^m\Bigr)u(x,0).
\la{gsol}
\ee
Since for constant coefficients $[\partial_x^n,\partial_x^m]=0$, the solution factorizes to 
\be
u(x,\dt)=\prod_{m=1}^M\e^{\dt a_m\partial_x^m}u(x,0),
\la{gsolf}
\ee
it is sufficient to study the effect of a single derivative operator at a time:
\be
u(x,\dt)=\e^{\dt a_m\partial_x^m}u(x,0).
\la{gsols}
\ee
Once numerical methods for solving (\ref{gsols}) for any $m$ are known,
the general equation (\ref{geq}) can be solved by a sequentual application of such schemes
according to (\ref{gsolf}). Higher dimension algorithms then follow from dimensional
splittings.

Conventionally, one expands out the RHS of (\ref{gsols}) and 
again approximates the spatial derivatives by finite-differences\cite{lw60,strang63}. 
However, in such an approach, how each derivative is to 
be approximated by a finite-difference (to what order, use which grid points) 
remained arbitrary and must be decided by some extrinsic considerations. 
Moreover, the resulting collection of schemes are just as disjoint and unrelated. 

This work proposes a theory of deriving explicit finite-difference schemes
that is still based on approximating the operator solution, but abandons the
practice of approximating derivatives {\it directly} by finite-differences.
Instead, such approximations are automatically generated by matching 
the finite-difference scheme to the formal solution and are completely prescribed 
by the temporal order of the algorithm. From this theory, 
all explicit finite-difference schemes for solving (\ref{gsols}) 
to any time-marching order are given by a single formula.  
 
The key idea of this theory is to use an operator form 
of the finite-difference scheme so that it can be transparently matched to
operator solution. This is described in the next section.
Once this this done, three fundamental theorems immediately follow that 
completely characterize all $n$th-order explicit time-marching algorithms for solving
any $m$-order partial differential equation.
In Section 3, the explicit form and the stability of the first-order time-marching algorithm 
are deterimined for all $m$ simultaneously. In Section 4, many higher-order time-marching algorithms 
are given for $m=1$ and $m=2$. These examples serve to illustrate the three theorems in Section 2. 
The nonlinear advection case is described in Section 5.
Some concluding remarks are given in Section 6. 

\section {Operator matching for linear equations}

An explicit finite-difference scheme seeks to approximate the exact solution (\ref{gsols}) via 
\be
u(x,\dt)=\sum_{i=1}^N c_i u(x+k_i\dx,0),
\la{fd}
\ee
where $\{k_i\}$ is a set of $N$ integers clustered around zero 
and $\{c_i\}$ is a set of coefficients. In the conventional approach,
$\{k_i\}$ and $\{c_i\}$ are by-products of the way spatial and temporal
derivatives are approximated by finite-differences,
and are therefore obtained concomitantly, mixed together.
This obscures the underlying relationship among 
schemes of different time-marching order $n$ for solving equations 
of different derivative order $m$. In this work, we disentangle the 
two and determine $\{k_i\}$ and $\{c_i\}$ separately.  

First, we will assume that $\{k_i\}$ is given set of $N$ integers, 
usually, a set of $N$ consecutive integers containing zero. The power of
our theory is that they need not be specified initially.
They are a set of parameters that will ultimately be
decided by the stability of the resulting numerical scheme.
 
Next, to determine $\{c_i\}$ for a given set of $\{k_i\}$, 
we make the following {\it key} observation:	  
that each grid value can also be represented in an operator form:  
\be
u(x+k_i\dx,0)=\e^{k_i\dx\partial_x}u(x,0).
\la{grid}
\ee
The finite-difference approximation (\ref{fd}) then corresponds to
\be
u(x,\dt)=\sum_{i=1}^N c_i\e^{k_i\dx\partial_x}u(x,0).
\la{fdap}
\ee
Comparing this to the exact solution (\ref{gsols}), due to the
linearity of the equation, the coefficients 
$c_i$ can be determined by solving the operator equality 
\be
\e^{\dt a_m\partial_x^m}=\sum_{i=1}^N c_i\e^{k_i\dx\partial_x}.
\la{opeq}
\ee
The simplest way to solve for $c_i$ is to Taylor expand both 
sides of (\ref{opeq}) and matches the powers of the derivative operator
$\partial_x$:
\be
1+\dt a_m \partial_x^m+\frac12 (\dt a_m)^2\partial_x^{2m}+\cdots
=\sum_{i=1}^N  c_i+\sum_{i=1}^N c_i(k_i\dx)\partial_x
+\frac12 \sum_{i=1}^N c_i(k_i\dx)^2\partial_x^2+\cdots .
\la{zeq}
\ee
From this, one immediately sees that for a given $m$,
an $n$th-order time-marching algorithm on the right, must
match up to the $nm^{th}$ power of $\partial_x$ on the left.
Thus $\{c_i\}$ must satisfy $N=nm+1$ linear order-conditions
and therefore requires the same number of grid points.
Thus we have proved the following fundamental theorem for
explicit finite-difference schemes:

\bt[Fundamental]
An $n$th-order time-marching finite-difference scheme of the form
$$
u(x,\dt)=\sum_{i=1}^N c_i u(x+k_i\dx,0)
$$
for solving the equation
$$
\frac{\partial u}{\partial t}=a_m\partial_x^m u,
$$
where $a_m$ is a real constant and $m$ a whole number $\ge 1$,
must have a minimum of $N=nm+1$ grid points. The latter is 
any set of $N$ integers $\{k_i\}$ clustered around zero.
\la{thefun}
\et

Note that this result is obtained without any prior knowledge of 
how derivatives are to be approximated by finite-differences.
Such approximations are automatically generated by the order conditions in 
(\ref{zeq}) for any set of $\{k_i\}$. In this theory of explicit finite-difference
schemes, everything follow from this set of order-conditions.   
One can easily check that all known low-order explicit schemes obey this
theorem. 

The set of order-conditions in (\ref{zeq}) can be solved easily,
and we have our central result:

\bt[Central]
An $n$th-order time-marching scheme
$$
u(x,\dt)=\sum_{i=1}^{N} c_i u(x+k_i\dx,0)
$$
with $\{c_i\}$ satisfying $N=nm+1$ order-conditions in (\ref{zeq})
for solving the equation
$$
\frac{\partial u}{\partial t}=a_m\partial_x^m u,
$$
has the closed-form solution
\be
c_i=\sum_{j=0}^n \frac{\nu_m^j}{j!}	L_i^{(jm)}(0)
\la{soluci}
\ee
where 
\be
\nu_m=\frac{\dt a_m}{\dx^m}
\la{cou}
\ee
is the generalized Courant number and $L_i^{(jm)}(0)$ 
are the $(jm)^{th}$-order derivatives of Lagrange polynominals of degree $N-1=nm$ 
$$
L_i(x)=\prod_{j=1 (\ne i)}^{{N}} \frac{(x-k_j)}{(k_i-k_j)},
\la{li}
$$
evaluated at the origin. 
\la{thecen}
\et
\bp
The order-condition (\ref{zeq}) reads individually,  
\ba
\sum_{i=1}^{N} c_i  &=& 1\la{seqg1}\\  
\sum_{i=1}^{N} c_ik_i^m  &=&\ m!\,\nu_m \nn\\ 
\sum_{i=1}^{N} c_ik_i^{2m}&=&(2m)!\,\frac{\nu_m^2}{2!}\nn\\ 
\cdots && \nn\\ 
\sum_{i=1}^{N} c_ik_i^{nm}&=&(nm)!\,\frac{\nu_m^n}{n!},
\la{seqg}
\ea
where $\nu_m$ is the generalized Courant number (\ref{cou}) and where 
all other powers of $k_i$ less than $nm$ sum to zero. We compare these
order conditions to  
the Vandermonde equation satisfied by Lagrange polynomials of $N-1=nm$ degree 
with grid points $\{k_i\}$ (See Appendix):
\ba
\begin{pmatrix}
 1     & 1        & 1      & ... & 1      \\
 k_1        & k_2      & k_3 & ...      & k_N \\
 k_1^2      & k_2^2      & k_3^2 & ...      & k_N^2 \\
 k_1^3      & k_2^3      & k_3^3 & ...      & k_N^3 \\
      ...      & ...           & ...      & ...      & ...  \\
 k_1^{(N-1)} & k_2^{(N-1)} & k_3^{(N-1)} & ... & k_N^{(N-1)} 
\end{pmatrix}
\begin{pmatrix}
L_1(x) \\ L_2(x)\\ L_3(x)\\ L_4(x)\\ ...\\ L_N(x)
\end{pmatrix}
=
\begin{pmatrix}
1\cr x\\ x^2\\ x^3\\ ...\\ x^{N-1} 
\end{pmatrix}
\la{solut2}
\ea
If we take the $\ell$th derivatives (including the zero-derivative) 
of this system of equations with respect to $x$ and set $x=0$ afterward, we would have
\be
\sum_{i=1}^N k_i^j L_i^{(\ell)}(0)=\ell!\delta_{j,\ell}
\la{leq}
\ee
This means that when $L_i^{(\ell)}(0)$ is summed over all powers 
of $k_i$ from 0 to $nm$, only the sum with $k_i^\ell$ is non-vanishing.
If $c_i$ were a sum of $L_i^{(\ell)}(0)$ terms, with $\ell=0,m,2m,\cdots, nm$, then
when $c_i$ is summed over powers of $k_i$, the sum will be
non-vanishing only at the required order-conditions 
(\ref{seqg1})-(\ref{seqg}). Adjusting the coefficients of $L_i^{(\ell)}(0)$ to
exactly match the order conditions then yields the solution (\ref{soluci}).
\ep 

Eq.(\ref{soluci}) is the master formula
for solving a general $m$-order partial differential equation to an arbitrary $n$th
time-marching order. All explicit schemes are related by their use of 
Lagrange polynomials. 
$L_i^{(\ell)}(0)$ is just $\ell!$ times the coefficient of the monomial
$x^\ell$ in $L_i(x)$.

For the next theorem, we need the
sums of $L_i(x)$ over $k_i^N$ and $k_i^{N+1}$, which are outside of (\ref{solut2}). They
are given by
\be
\sum_{i=1}^Nk_i^N L_i(x)=P(x) \qquad{\rm and }\qquad \sum_{i=1}^Nk_i^{N+1} L_i(x)=Q(x),
\ee
where 
\be
P(x)=x^N-s(x),\quad	Q(x)=x^{N+1}-\Bigl(x+\sum_{i=1}^Nk_i\Bigr)s(x),
\quad{\rm and }\quad s(x)=\prod_{i=1}^N(x-k_i).
\ee
Note that $P(x)$ and $Q(x)$ are just $N-1$ degree polynomials.
Taking the $\ell^{th}$ derivative ($0\le\ell\le N-1$) with respect to $x$ and set $x=0$ afterward yields,
\be
\sum_{i=1}^Nk_i^N L_i^{(\ell)}(0)=P^{(\ell)}(0)
\qquad{\rm and}\qquad \sum_{i=1}^Nk_i^{N+1} L_i^{(\ell)}(0)=Q^{(\ell)}(0)
\la{np}
\ee

By the way these schemes are constructed, it is very easy to compute their errors
with respect to the exact solution. Moreover, all such
explicit finite-difference schemes are characterized by a uniformity property:
\bt
An $n$th-order time-marching scheme
$$
u(x,\dt)=\sum_{i=1}^{N} c_i u(x+k_i\dx,0)
$$
with $N=nm+1$ and with $\{c_i\}$ given by Theorem \ref{thecen}, 
for solving the equation
$$
\frac{\partial u}{\partial t}=a_m\partial_x^m u,
$$
approximates all spatial derivatives $u^{(jm)}(x,0)$ from $j=0$ to $j=n$
uniformly to order $\dx^{nm}$ and has an overall local error of
\ba
&&E=\dx^{N}\Bigl(\sum_{j=0}^n \frac{\nu_m^j}{j!}P^{(jm)}(0)\Bigr) \frac{u^{(N)}(x,0)}{N!}
+\dx^{N+1}\Bigl(\sum_{j=0}^n \frac{\nu_m^j}{j!}Q^{(jm)}(0)\Bigr) \frac{u^{(N+1)}(x,0)}{(N+1)!}\nn\\
&&\qquad\qquad\qquad\qquad\qquad\qquad\qquad\qquad\qquad\qquad-(\dt a_m)^{n+1}\frac{u^{((n+1)m)}(x,0)}{(n+1)!}. 
\la{te}
\ea
The local truncation error is just $E/\dt$.
\la{theerr}
\et
\bp
Substitute in the solution for $c_i$ from (\ref{soluci}) gives
\ba
&&\sum_{i=1}^{N} c_i u(x+k_i\dx,0)
        =\sum_{j=0}^n \frac{\nu_m^j}{j!}
            \sum_{i=1}^{N} L_i^{(jm)}(0)u(x+k_i\dx,0)\nn\\
  &&\quad=\sum_{j=0}^n \frac{\nu_m^j}{j!}
    \sum_{i=1}^{N} L_i^{(jm)}(0)\Bigl[u(x,0)+(k_i\dx) u^{(1)}(x,0)+\frac1{2!}(k_i\dx)^2 u^{(2)}(x,0) +\cdots
	\nn\\
&&\qquad\qquad\qquad+\frac{(k_i\dx)^{N}}{N!}u^{(N)}(x,0)+\frac{(k_i\dx)^{N+1}}{(N+1)!} u^{(N+1)}(x,0)+O(\dx^{N+2})\Bigr]
\ea
By (\ref{leq}) and (\ref{np}), we have
\ba
&&\sum_{i=1}^{N} c_i u(x+k_i\dx,0)
=\sum_{j=0}^n \frac{\nu_m^j}{j!}
\Bigl[(\dx)^{jm} u^{(jm)}(x,0) +\frac{P^{(jm)}(0)}{N!} \dx^{N} u^{(N)}(x,0)\nn\\
&&\qquad\qquad\qquad\qquad\qquad\qquad\quad+\frac{Q^{(jm)}(0)}{(N+1)!} \dx^{N+1} u^{(N+1)}(x,0)+O(\dx^{N+2})\Bigr].
\la{unf}
\ea
All approximations of $u^{(jm)}(x,0)$ for $0\le j\le n$ are therefore uniformily correct to at least
spatial order $N-1=nm$. Subtracting the exact solution (\ref{gsols}) from above gives the local error (\ref{te}).
\ep 

Note that $P^{(jm)}(0)$ may vanish, if so, that derivative
approximation will then be correct to one order higher. 
This is why we needed the $Q^{(jm)}(0)$ term for the diffusion equation.
Also, the $j=0$ case means that $u(x,0)$ is correctly approximated to order $\dx^{nm}$,
if $\{k_i\}$ does not contain 0. See (\ref{lzero}) below. 

This theorem states that the error analysis can be done once for all explicit 
finite-differences schemes. There is no need to do Taylor expansions for 
each finite-difference scheme, one by one.
Also, this theorem shows that there is no arbitrariness in specifying the order of the 
spatial derivatives approximations. At time-marching order $n$, all spatial derivatives 
(including the initial function itself)
must be uniformly approximated to order $nm$.
The spatial order of approximation is completely fixed
by the temporal order of the algorithm. Examples 
illustrating these three theorems will be given in Section 4.

\section {Complete characterization of first-order algorithms}

From Theorem \ref{thecen}, all explicit numerical schemes
are given by the master formula (\ref{soluci}).
However, for a given $(m,n)$, it is easy to show that for $(1,n)$ and $(m,1)$,
the coefficients $c_i$ are particularly simple. 
We will discuss the first case in the next section. For the second case, 
the set of $\{c_i\}$ has the following simple form:

\bt
The $m+1$ first-order time marching finite-difference schemes 
$$
u(x,\dt)=\sum_{k=-r}^{m-r} c_k u(x+k\dx,0)
$$
characterized by $r=0,1,2,...m$	for solving the equation
$$
\frac{\partial u}{\partial t}=a_m\partial_x^m u
$$
have explicit solutions
\be
c_0=1+(-1)^m(-1)^{r} C_{r}^m \nu_m\quad {\rm and}\quad
c_k=(-1)^m(-1)^{r+k} C_{r+k}^m \nu_m
\la{solu}
\ee
given in terms of the generalized Courant number $\nu_m=\dt a_m/\dx^m$
and binomial coefficients
\be
C_k^m=\frac{m!}{k!(m-k)!}.
\ee
\et
\bp
For a first-order time-marching scheme, we have 
$N=m+1$ grid points, which we can take to be 
$k_i=\{-r,-r+1,\dots -1,0,1,\dots s\}$, where $s=m-r$
and where each value of $r=0,1,2,...m$ labels a distinct algorithm.
The corresponding coefficients can then be denoted directly by their
$k_i$ values as 
$c_{-r},c_{-r+1}, \dots, c_0, c_1, \dots c_s $. From Theorem 3,
since each Lagrangian polynomial is defined by
\be
L_i(x)=\prod_{j=1 (\ne i)}^{m+1} \frac{(x-k_j)}{(k_i-k_j)},
\ee
one has
\be
L_i^{(m)}(0)=\frac{m!}{\prod_{j=1 (\ne i)}^{m+1} (k_i-k_j)}.
\ee
We now eliminate the index ``$i$" by replacing $k_i$
by its actual value denoted by $k$. 
One then sees that
\be
L_k(0)=\prod_{j=-r (\ne k)}^{s} \frac{(0-j)}{(k-j)}=\delta_{k,0}
\la{lzero}
\ee
and
\ba
L_{k}^{(m)}(0)&=&\frac{(-1)^m m!}{(s-k)(s-k-1)...(-r-k+1)(-r-k)}\nn\\
&=&\frac{(-1)^m m!}{(s-k)! (-1)...(-r-k+1)(-r-k)}\nn\\
&=&\frac{(-1)^m(-1)^{r+k} m!}{(s-k)!(r+k)!}
=\frac{(-1)^m(-1)^{r+k}m!}{(m-r-k)!(r+k)!}\nn\\
&=&(-1)^m(-1)^{r+k}C^m_{r+k},
\ea 
which produces the explicit solution (\ref{solu}).
\ep

To gain insights about this set of first-order algorithms for all $m$, consider
the generation function for the coefficients $c_k$:
\ba
g(x)&=&\sum_{k=-r}^{m-r} c_kx^k=1+\nu_m \sum_{k=-r}^{m-r} (-1)^m(-1)^{r+k} C_{r+k}^m x^k\nn\\
&=&1+\nu_m x^{-r}\sum_{k=-r}^{m-r} (-1)^m(-1)^{r+k} C_{r+k}^m x^{r+k}
\ea 
Shifting the dummy variable $k\rightarrow r+k$ gives
\ba
g(x)&=1+&\nu_m x^{-r}\sum_{k=0}^m (-1)^m(-1)^{k} C_{k}^m x^{k}\nn\\
&=1+&\nu_m x^{-r}(x-1)^m.
\la{genf}
\ea
Thus the coefficients of the algorithm are just coefficients of $(x-1)^m$. 

This generation function can now be used to determine the stability 
of this set of first-order algorithms for all $m$ simultaneously via the following
two theorems.

\bt
If the first-order time marching finite-difference scheme 
described in Theorem 4 for solving the equation
$$
\frac{\partial u}{\partial t}=a_m\partial_x^m u
$$
is von-Neumman stable, then its range of stability is limited to
$$
|\nu_m|\le \frac1{2^{m-1}}.
$$
\et
\bp
The generation function (\ref{genf}) give the following amplification factor 
for a single Fourier mode $\e^{ipx}$,
\be
g=\sum_{k=-r}^{m-r} c_k (\e^{i\theta})^k=1+\nu_m \e^{-i r \theta}(\e^{i\theta}-1)^m,
\ee
where we have denoted $\theta=p\dx$. since
\be
\e^{i\theta}-1=\e^{i\theta/2}2i\sin(\theta/2)=\e^{i(\theta/2+\pi/2)}2\sin(\theta/2)
\ee
we have
\be
g=1+\nu_m \e^{i(m(\theta+\pi)/2-r\theta)}[2\sin(\theta/2)]^m,
\ee
and therefore
\be
|g|^2=1+2\cos(\Phi)\nu_m[2\sin(\theta/2)]^m+\nu_m^2[2\sin(\theta/2)]^{2m}
\ee
with 
\be
\Phi=\frac\theta{2}(m-2r)+m\frac{\pi}2.
\ee 						
The algorithm can be stable at small $|\nu_m|$ only if 
\be
{\rm sgn}(\nu_m)\cos(\Phi)< 0
\la{stabc}
\ee 
for all $\theta\in [0,2\pi]$. In this case, $|g|^2$ as a quadratic function of $|\nu_m|$
would first dip below one, reaching a minimum at $|\nu_m|_{min}=|\cos(\Phi)|/[2\sin(\theta/2)]^m$,
then backs up to one at $2|\nu_m|_{min}$. Thus the stability range of $|\nu_m|$ is limited
to
\ba
 |\nu_m|&& \le 2|\nu_m|_{min}=\frac{2|\cos(\Phi)|}{[2\sin(\theta/2)]^m} \nn\\
         && \le \frac1{2^{m-1}},
\ea
since the growth of $|g|^2$ is the greatest along $\theta=\pi$ with
$|\cos(\Phi)|=|\cos((m-r)\pi)|=|\cos(s\pi)|=1$.
\ep

Theorem 5 ``explains" why the upwind algorithm for solving the
$m=1$ advection equation is stable only for $|\nu_1|\le 1$
and that the $m=2$ diffusion algorithm is stable only for
$|\nu_2|\le 1/2$. These are not just isolated idiosyncrasies of
individual algorithm; they are part of the pattern 
of stability mandated by Theorem 5. One can easily check that this
theorem is true for other values of $m$. Thus with increasing $m$, the
range of stability decreases geometrically. 

We can now decide, 
among the $m+1$ first-order algorithms corresponding to $r=0,1,2,\dots m$ of
Theorem 4, which one is von-Neumman stable. Surprisingly,  
there is {\it at most} one stable first-order algorithm for a given value of 
$m$ and the sign of $\nu_m$: 
 
\bt
Among the $m+1$ first-order time marching finite-difference schemes 
described in Theorem 4 for solving the equation
\be
\frac{\partial u}{\partial t}=a_m\partial_x^m u,
\la{meq}
\ee
there is at most one stable algorithm for each value of $m$ and the sign of $a_m$.
For $m=2\ell$ the algorithm $r=\ell$ is stable only for {\rm sgn}$(a_m)=(-1)^{\ell-1}$.
For {\rm sgn}$(a_m)=(-1)^{\ell}$, there are no stable algorithms.
For $m=2\ell-1$, the algorithms $r=\ell$ and $r=\ell-1$ are stable 
for {\rm sgn}$(a_m)=(-1)^\ell$ and {\rm sgn}$(a_m)=(-1)^{\ell-1}$ respectively.
\et
\bp
Consider first the even case of $m=2\ell$, with $\ell=1,2,3,\cdots$. In this
case
\be
\cos(\Phi)=\cos(\theta(\ell-r)+\ell\pi)=(-1)^\ell\cos(\theta(\ell-r)).
\ee
If $(\ell-r)\ne0$, then as $\theta$ ranges from 0 to $2\pi$, $\cos(\theta(\ell-r))$
must change sign and the stability condition (\ref{stabc}) cannot hold for
all values of $\theta$. The only possible stable algorithm is therefore the central-symmetric
algorithm with $r=\ell$, which then places the following restriction on the sign of $\nu_m$:
\be
{\rm sgn}(\nu_m)(-1)^\ell=-1 \quad\Longrightarrow\quad {\rm sgn}(\nu_m)=(-1)^{m/2-1}.
\ee
That is, a stable first-order algorithm is only possibe for $a_2>0$, $a_4<0$, $a_6>0$, etc.,
and no stable algorithm otherwise.

For the odd case of $m=2\ell-1$, with $\ell=1,2,3,\cdots$, we now have
\be
\cos(\Phi)=\cos\Bigl[\theta(\ell-r)+\ell\pi-\frac12 (\theta+\pi) \Bigr]
=(-1)^\ell\cos\Bigl[\theta(\ell-r)-\frac12 (\theta+\pi)\Bigr].
\ee 
For $r=\ell$, $\cos\bigl[-\frac12 (\theta+\pi)\bigr]<0$ for $0<\theta<2\pi$. Thus this algorithm
is stable for
\be
{\rm sgn}(\nu_m)(-1)^{\ell+1}=-1 \quad\Longrightarrow\quad {\rm sgn}(\nu_m)=(-1)^{\ell}.
\ee
For $r=\ell-1$, we have $\cos\bigl[\frac12 (\theta-\pi)\bigr]>0$ for $0<\theta<2\pi$ and this
algorithm is stable for
\be
{\rm sgn}(\nu_m)(-1)^{\ell}=-1 \quad\Longrightarrow\quad {\rm sgn}(\nu_m)=(-1)^{\ell-1}.
\ee
Other than these two values of $r$, 
$\cos\bigl[\theta(\ell-r)-\frac12 (\theta+\pi)\bigr]$ will changes sign as $\theta$ ranges over $2\pi$.
For each sign of $a_{2\ell-1}$, there is only one stable algorithm.
\ep

Theorems 1, 2, 4, 5, 6 completely characterize all minimum grid-point, first-order time-marching
algorithms for solving (\ref{meq}). The pattern of stability proscribed by Theorem 6
is easily understood from the following plane wave solution to (\ref{meq}):
\be
u(x,t)=A\e^{ipx+a_m(ip)^m t}.
\ee
For $m=1$, the wave propagates to the positive x-direction for $a_1<0$, therefore
only the algorithm $r=1$ is stable, corresponding to the upwind algorithm. For
$a_1>0$, the wave propagates to the negative x-direction and $r=0$ is the corresponding
upwind algorithm. 

For $m=2$, the wave decays in time only for $a_2>0$ and $r=1$ gives the well-known
first-order diffusion algorithm. For $a_2<0$, there is no stable algorithm 
because the solution grows without bound with time. 

For $m=3$, the wave propagates to the positive x-direction with $a_3>0$; therefore
only the $r=2$ algorithm is stable, the analog of the $m=1$ upwind algorithm. For
$a_3<0$, the wave propagates to the negative x-direction and $r=1$ is the analogous
upwind algorithm. Note that the scheme with 
grid points $\{k_i\}=\{-2,-1,1,2\}$, excluding $0$, is {\it unstable}. 

For $m=4$, the wave decays in time only for $a_4<0$ and pattern repeats as $i^m$
cycles throught its four possible values.

\section {Higher order time-marching algorithms}
   
For the case of $(m,n)=(1,n)$, the algorithm for solving the advection equation 
to the $n$th time-marching order is given by 									  
\be
c_i=\sum_{k=0}^n L_i^{(k)}(0)\frac{\nu_1^k}{k!}=L_i(\nu_1),
\la{anyad}
\ee
which is the seminal case studied by Strang\cite{strang62}, 
Iserles and Strang\cite{is83} and recounted in Ref.\cite{hv}.
This is now just a special case of Theorem \ref{thecen}. Note that 
the last equality can be used to identify $L_i^{(k)}(0)$ needed 
by other algorithms, see further discussion below.   
 
For $k_i=\{-r,-r+1,\dots -1,0,1,\dots s\}$, as shown by Strang\cite{strang62},
and Iserles and Strang\cite{is83}, 
only three cases are stable for each sign of $a_1$.
For the conventional choice of $a_1<0$, where the wave propagates from left to right,
 $r=s+1$ and $r=s$ are stable for $0\le |\nu_1|\le 1$, and $r=s+2$ is stable for $0\le |\nu_1|\le 2$.
Since $r+s=n$, the order of each type of algorithm are $n=2s+1$, $n=2s$ and $n=2s+2$ respectively.
Thus there is one stable algorithm at each odd-order and two stable algorithms at each even-order.
The odd-order algorithms begin with the first-order upwind (UW) scheme with $s=0$
and the even-order schemes begin with the second-order Lax-Wendroff\cite{lw60}(LW) 
scheme with $s=1$ and the second-order Beam-Warming\cite{bw76} (BW) 
scheme with $s=0$, respectively. The even-order schemes can be distinguished
as being of the LW-type ($n=2s$) or BW-type ($n=2s+2$). 

The local error in this case is particularly simple. From (\ref{te}),
for $m=1$ and $N=n+1$, we have (ignoring the $Q^{jm}(0)$ term),
\be
E=\Bigl[\dx^{n+1}\Bigl(\sum_{j=0}^n\frac{(\nu_1)^j}{j!}P^{(j)}(0)\Bigr)-(\dt a_1)^{n+1}\Bigr] 
\frac{u^{(n+1)}(x,0)}{(n+1)!}
\la{tea}
\ee
Since $P(x)$ is a polynomial of $N-1=n$ degree, we have
\ba
E
&=&\Bigl[\dx^{n+1}P(\nu_1)-(\dt a_1)^{n+1}\Bigr] \frac{u^{(n+1)}(x,0)}{(n+1)!}\nn\\
&=&\Bigl[\dx^{n+1}\Bigl(\nu_1^{n+1}-\prod_{i=1}^{n+1}(\nu_1-k_i)\Bigr)-(\dt a_1)^{n+1}\Bigr] 
\frac{u^{(n+1)}(x,0)}{(n+1)!}\nn\\
&=&-\dx^{n+1}\prod_{i=1}^{n+1}(\nu_1-k_i)\frac{u^{(n+1)}(x,0)}{(n+1)!}.
\la{teb}
\ea
The local truncation error is obtained by dividing the above by $\dt$.

For later illustration purposes, we list below the third-order scheme 
corresponding to $s=1$ with
$\{k_i\}=\{-2,-1,0,1\}$. The coefficients directly from (\ref{anyad}) are
\be
c_{-2}=-\frac{\nu_1}{6}(\nu_1^2-1)\qquad c_{-1}=\frac{\nu_1}{2}(\nu_1+2)(\nu_1-1)
\ee 
\be
c_{0}=-\frac12(\nu_1+2)(\nu_1^2-1) \qquad c_{1}=\frac{\nu_1}{6}(\nu_1+2)(\nu_1+1),
\ee 
and the algorithm can also be arranged as a sum over powers of $\nu_1$:
\ba
u^{k+1}_j
=u^k_j&+&\nu_1(\frac1{6}u^k_{j-2}-u^k_{j-1}+\frac12u^k_{j}+\frac13 u^k_{j+1})\nn\\
&+&\frac{\nu_1^2}{2!}(u^k_{j-1}-2u^k_{j}+u^k_{j+1})\nn\\
&+&\frac{\nu_1^3}{3!}(-u^k_{j-2}+3u^k_{j-1}-3u^k_{j}+ u^k_{j+1}),
\la{ad3}
\ea
where we have denoted $u^k_j\equiv u(k\dt,j\dx)$.
Similarly, the fourth-order LW-type algorithm corresponding to $\{k_i\}=\{-2,-1,0,1,2\}$
can be arranged as
\ba
u^{k+1}_j
=u^k_j&+&\nu_1(\frac1{12}u^k_{j-2}-\frac23 u^k_{j-1}+\frac23 u^k_{j+1}-\frac1{12} u^k_{j+2})\nn\\
&+&\frac{\nu_1^2}{2!}(-\frac1{12} u^k_{j-2}+\frac43 u^k_{j-1}-\frac52 u^k_{j}+\frac43 u^k_{j+1}-\frac1{12} u^k_{j+2})\nn\\
&+&\frac{\nu_1^3}{3!}(-\frac12 u^k_{j-2}+u^k_{j-1}-u^k_{j+1}+\frac12 u^k_{j+2})\nn\\
&+&\frac{\nu_1^4}{4!}(u^k_{j-2}-4u^k_{j-1}+6u^k_{j}-4 u^k_{j+1}+u^k_{j+2}).
\la{ad4}
\ea
Each parenthese in (\ref{ad3}) and (\ref{ad4}) gives the corresponding
third and fourth-order spatial discretization
of derivatives respectively, as mandated by Theorem 2. The coefficients of the last term 
are just those of $(x-1)^3$ and $(x-1)^4$, in accordance with Theorem 4.

The above examples are for later illustrations only. In practice, it is {\it absolutely unnecessary} to write 
out the coefficients $c_i$ explicitly as in the above examples, or disentangle them into
powers of $\nu_1$. It is only necessary to write a short routine to compute $c_i$ directly
from (\ref{anyad}) for a given set of $\{k_i\}$ and generate an algorithm of any order.
This is illustrated below. 

\begin{figure}[hbt]
\includegraphics[width=0.49\linewidth]{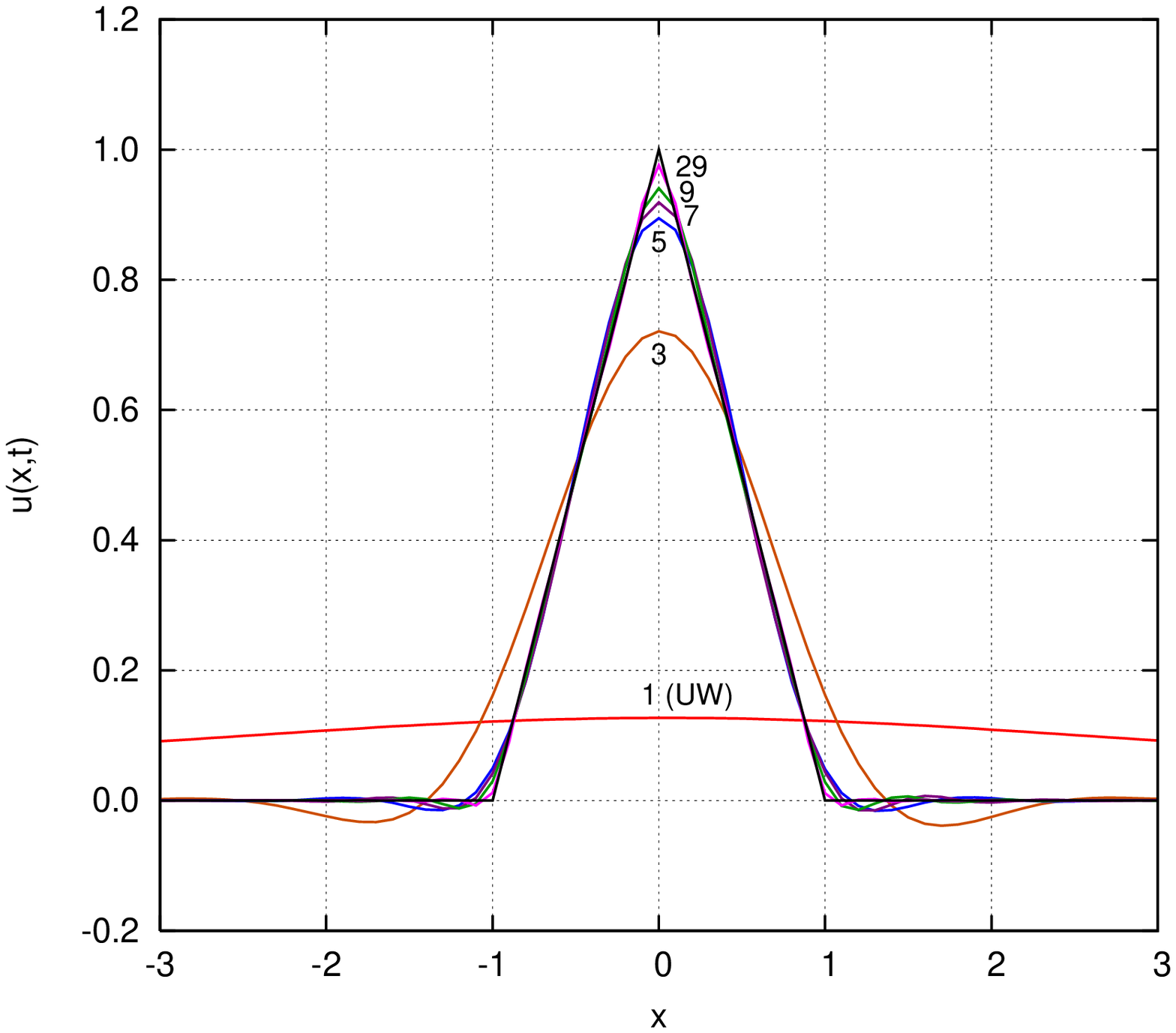}
\includegraphics[width=0.49\linewidth]{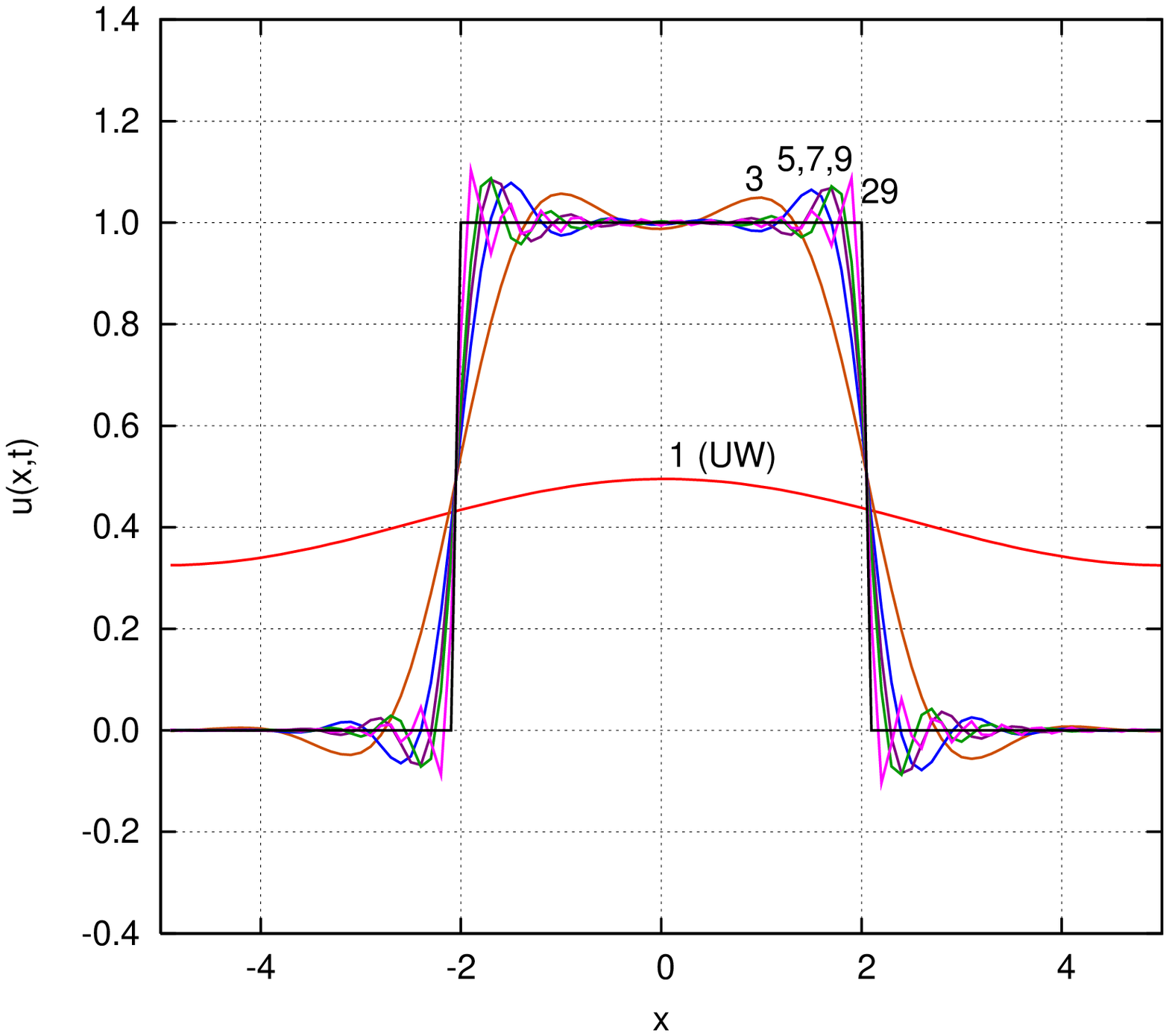}
\caption[]{\label{aduw} 
The propagation of an initial triangular and a rectangular profile (in black) 
50 times around a periodic box of [-5,5] with $\dx=0.1$, $\dt=0.08$, $v=1$, and $\nu_1=0.8$,
corresponding to 6250 iterations of each algorithm. The numbers label odd-order algorithms
beginning with the first-order upwind (UW) scheme. The highest order is 29.
}
\end{figure}

\begin{figure}[hbt]
\includegraphics[width=0.49\linewidth]{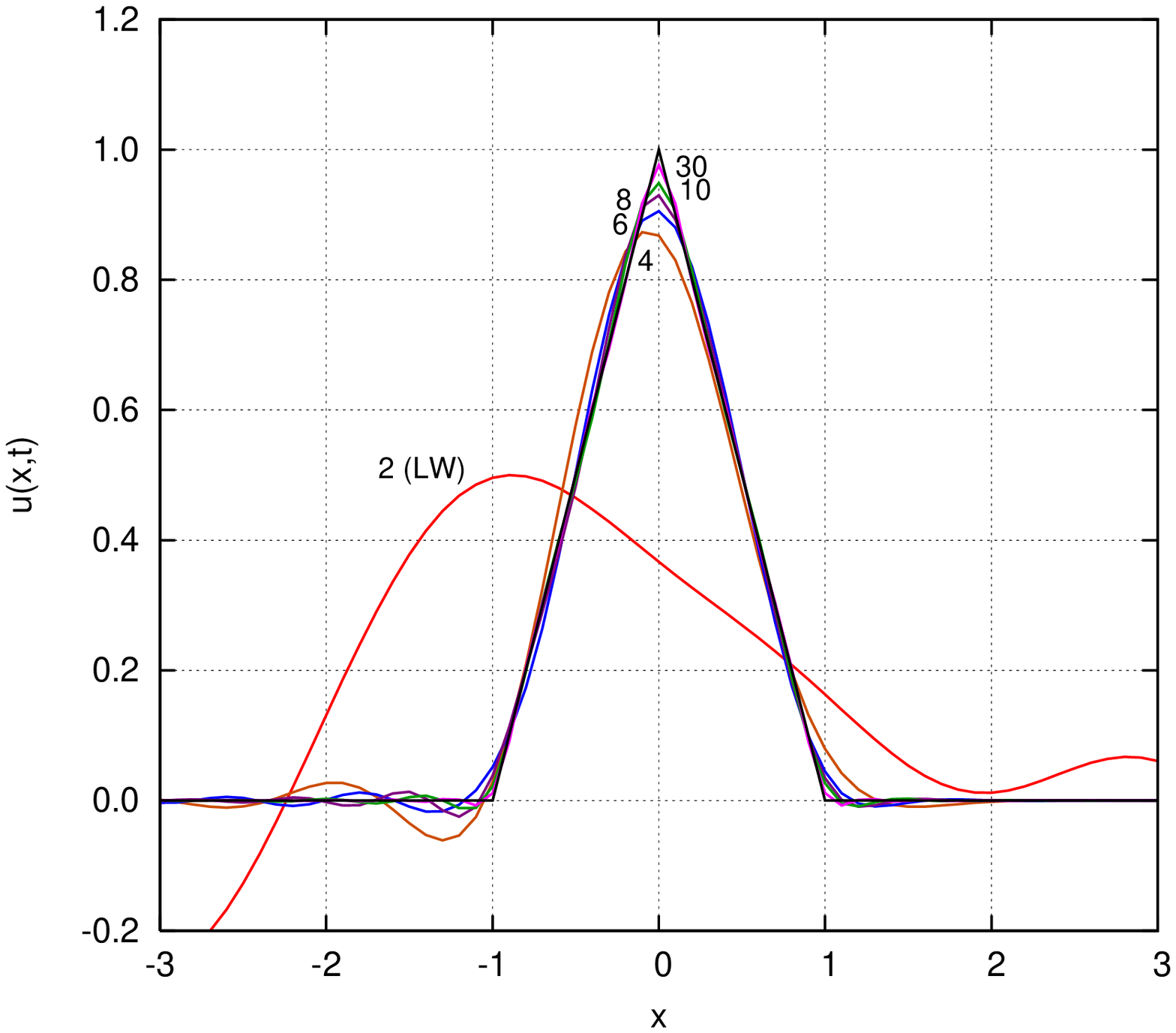}
\includegraphics[width=0.49\linewidth]{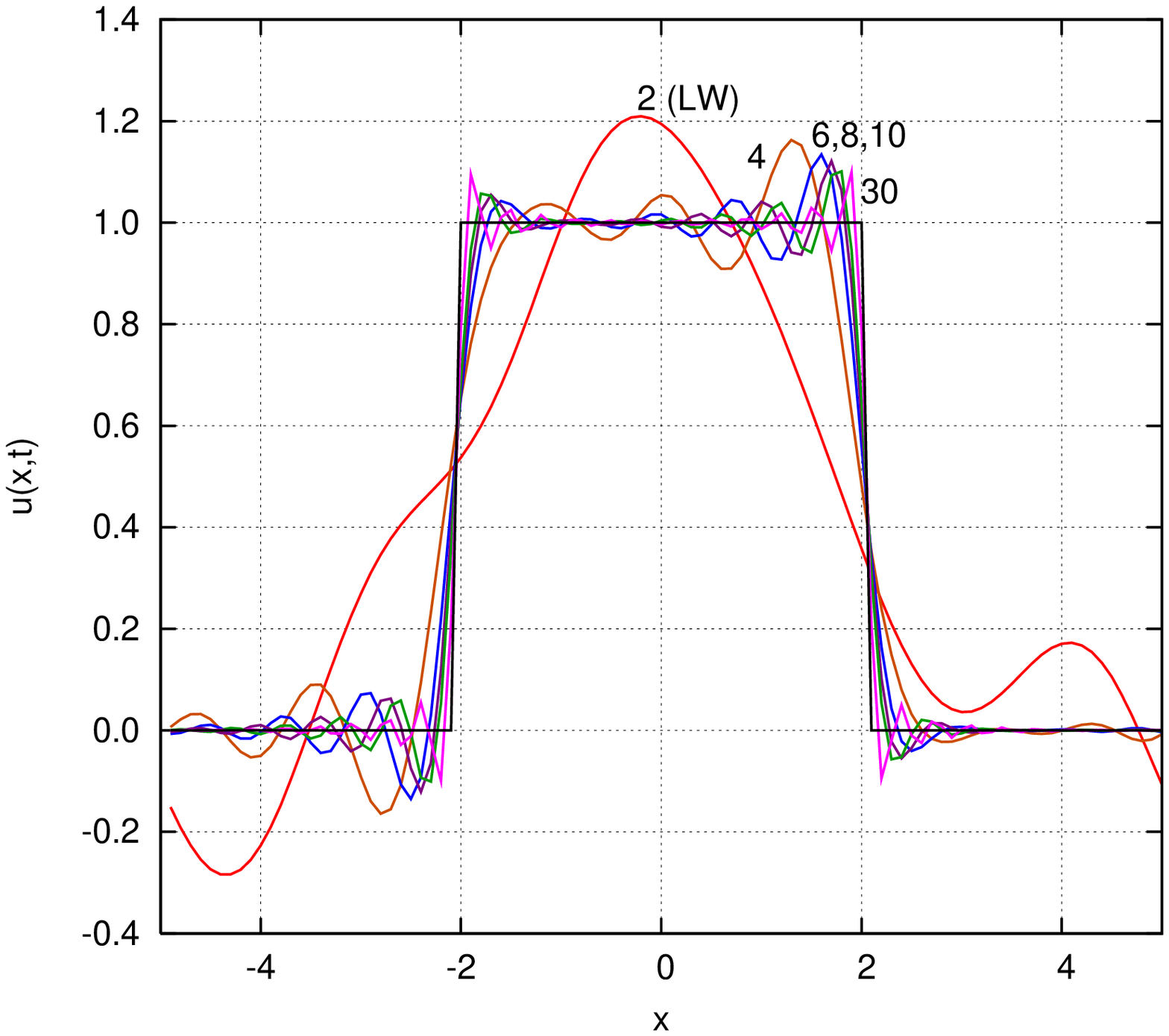}
\caption[]{\label{adlw} 
Same as Fig.\ref{aduw} for even-order algorithm whose lowest-order member is
the second-order Lax-Wendroff (LW) scheme. The highest order here is 30. 
}
\end{figure}

\begin{figure}[hbt]
\includegraphics[width=0.49\linewidth]{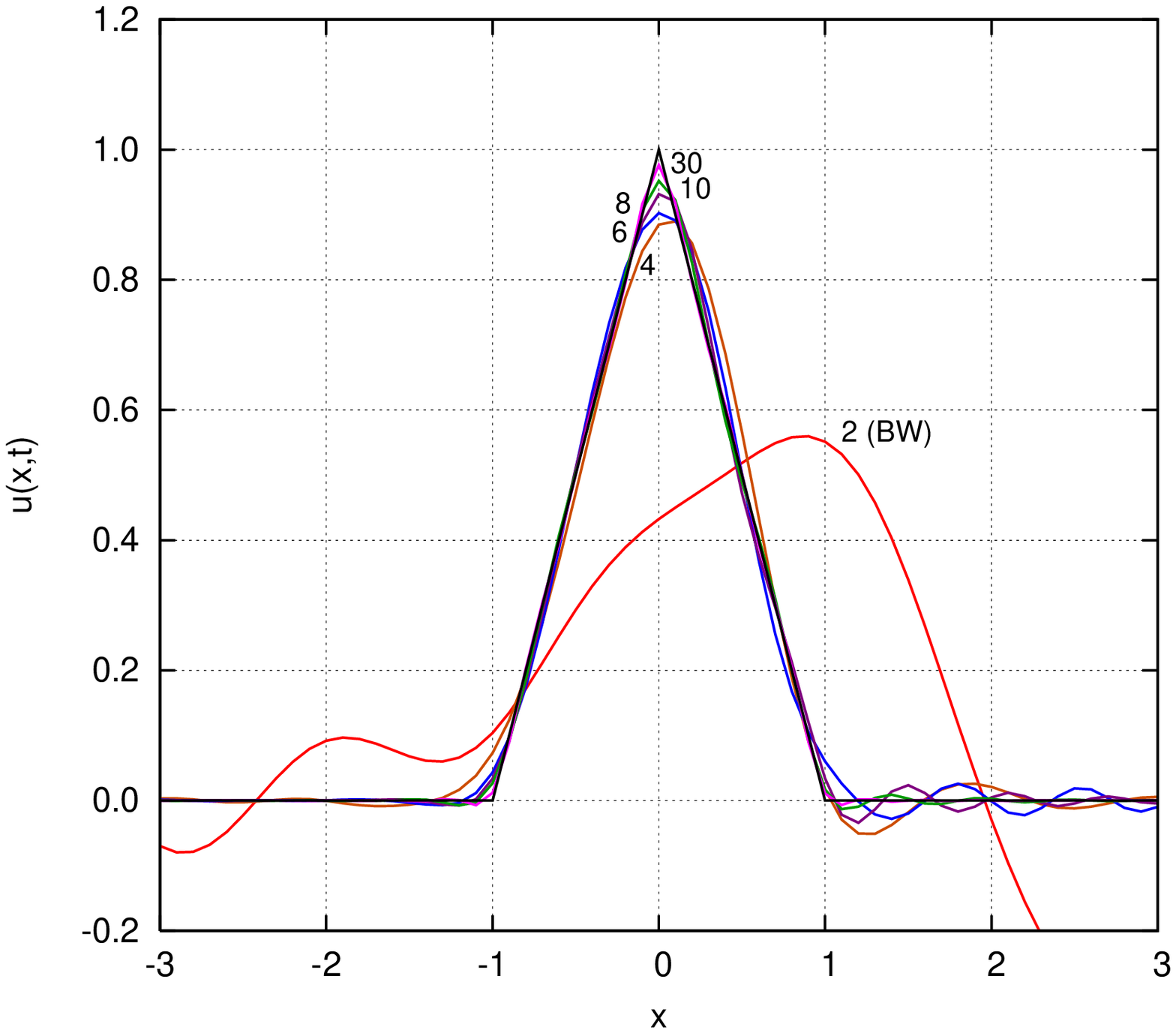}
\includegraphics[width=0.49\linewidth]{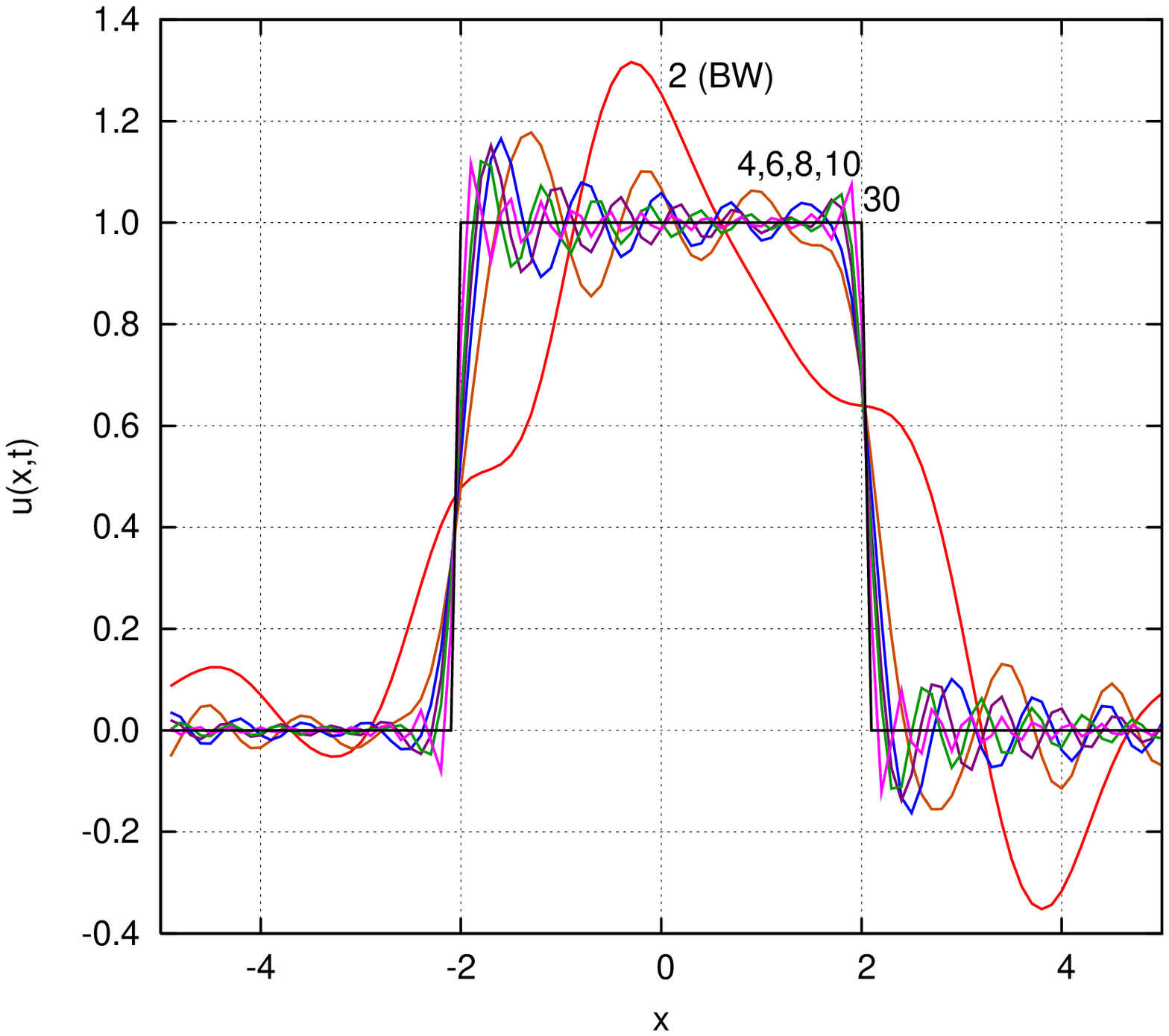}
\caption[]{\label{adbw} 
Same as Fig.\ref{aduw} for even-order algorithm whose lowest-order member is
the second-order Beam-Warming (BW) scheme. The highest order here is 30. 
}
\end{figure}

In solving the advection equation, it is well-known that low-order algorithms 
are plagued with unwanted damping and negative oscillations. In the following
figures, we examine how each type of algorithm converges toward the exact
solution with increasing time-marching orders. Figs.\ref{aduw}, \ref{adlw} and \ref{adbw}
show results for the odd-order UW-type and the two even-order LW-type and BW-type algorithms respectively. 
In each case, the lowest order algorithm, UW, LW and BW are all damped and dispersed beyond
recognition. At increasing order (at fixed $\dx$ and $\dt$), the convergence toward 
the undamped triangular profile is
excellent. For the rectangular profile, the convergence is consistent with
having Gibb's oscillations. The odd-order algorithms 
preserve the left-right symmetry of the original profile despite oscillations, 
whereas the two even-order algorithms are marred by asymmetries and phase errors 
until very high orders. In this study, the order 5 scheme seemed optimal.
Beyond order 5, the improvement is incremental. The 29th and 30th-order algorithms are of course
not very practical. They are shown here just to illustrate the fact that very high order 
algorithms are possible.

For solving the diffusion equation with $m=2$,
it is natural to take $N=2n+1$ grid points to be 
$\{k_i\}=\{-n,...-1,0,1,2..,n\}$ 
with $c_{-i}=c_i$. For $n=1$ one obtains the familiar 
first-order time-marching algorithm from Theorem 4:
\be
 u^{k+1}_j=u^k_j+\nu_2( u^k_{j-1}-2u^k_{j}+u^k_{j+1}).
\ee
The coefficients multiplying
$\nu_2$ are just those of $(x-1)^2$.
For time-marching order 2, one has
\be
c_0=1-\frac52 \nu_2+3 \nu_2^2\qquad c_1=\frac43 \nu_2-2\nu_2^2\qquad c_2=-\frac1{12}\nu_2
+\frac12 \nu_2^2
\la{diff2}
\ee
and the resulting algorithm is
\ba
 u^{k+1}_j=u^k_j
 &+&\nu_2(-\frac1{12} u^k_{j-2}+\frac43 u^k_{j-1}-\frac52 u^k_{j}+\frac43 u^k_{j+1}-\frac1{12} u^k_{j+2})\nn\\
&+&\frac{\nu_2^2}{2!}(u^k_{j-2}-4u^k_{j-1}+6u^k_{j}-4 u^k_{j+1}+u^k_{j+2}).
 \la{secd}
\ea
Comparing this to the fourth-order advection algorithm (\ref{ad4}), one sees that the coefficients
inside the parentheses are just those of the second and fourth order terms of (\ref{ad4}). 
Thus the coefficients of the diffusion algorithm are simply those of the even order terms of the
advection algorithm, with appropriate change of factors $\nu_1^{2k}/(2k)!\rightarrow\nu_2^k/k!$,
provided that both are using the same set of $\{k_i\}$. (Similarly, one can pick out the 
$jm$-order terms of the advection scheme to generate algorithms for solving the $m$-order equation.)   

For order 3, one has
\be
c_0=1-\frac{49}{18} \nu_2+\frac{14}3 \nu_2^2 -\frac{10}3 \nu_2^3  \qquad 
c_1=\frac{3}{2} \nu_2-\frac{13}4 \nu_2^2 +\frac{5}2 \nu_2^3
\ee
\be
c_2=-\frac{3}{20} \nu_2+\nu_2^2 -\nu_2^3\qquad\qquad\quad
c_3=\frac{1}{90} \nu_2-\frac1{12} \nu_2^2 +\frac{1}6 \nu_2^3.
\ee
Again, the coefficients of $\nu_2^3$ are now
binominal coefficients of $(x-1)^6$ multiplied by $1/3!$. The coefficients of 
$\nu_2$ and $\nu_2^2$ are from 
Theorem 3.
The algorithm is now correct to $sixth-order$ in spatial discretizations.

For order 4, one has
\be
c_0=1-\frac{205}{72} \nu_2+\frac{91}{16} \nu_2^2 -\frac{25}4 \nu_2^3 +\frac{35}{12} \nu_2^4 
\ee
\be
c_1=\frac{8}{5} \nu_2-\frac{61}{15} \nu_2^2 +\frac{29}6 \nu_2^3	-\frac{7}3 \nu_2^4\qquad
c_2=-\frac{1}{5} \nu_2+\frac{169}{120}\nu_2^2 -\frac{13}{6}\nu_2^3+\frac{7}6 \nu_2^4
\ee
\be
c_3=\frac{8}{315} \nu_2-\frac1{5} \nu_2^2 +\frac{1}2 \nu_2^3-\frac{1}3 \nu_2^4\qquad
c_4=-\frac{1}{560} \nu_2+\frac{7}{480}\nu_2^2 -\frac{1}{24}\nu_2^3+\frac{1}{24} \nu_2^4,
\ee
which is correct to the {\it eighth-order} in spatial discretizations.

In all these cases, one can check that the local error is correctly given by (\ref{te}),
unfortunately, there does not seem to be a closed form for the sum over $Q^{(jm)}(0)$,
and hence no simple expression for the local error as in (\ref{teb}). 
  
For these algorithms, the amplification factor for a single Fourier mode $\e^{ipx}$ is
\be
g=1-4\sum_{j=1}^n c_j\sin^2(\frac{j}2 \theta)
\ee
where $\theta=p\dx$. For $n=1$, $c_1=\nu_2$, one obtains the usual stability criterion of 
$\nu_2\le \nu_c$,
where the critical stability point is $\nu_c=1/2$.
If one simply increases the spatial order of discretizing $\partial^2u^n_j$ to fourth-order
without increasing the temporal order, $\nu_c$ decreases to $3/8=0.375$. (This is
algorithm (\ref{secd}) without the second-order $\nu_2^2$ term.) 
However, the full second-order time-marching algorithm
(\ref{secd}) has {\it increased} stability, with $\nu_c=2/3=0.667$. 
Similarly, the third and fourth-order 
time-marching algorithm have increased stability of $\nu_c=0.841$ and $\nu_c=1.015$ respectively, 
while keeping only the first-order term in $\nu_2$ have decreased stability of 
$\nu_c=45/136=0.331$ and $\nu_c=315/1024=0.308$. Thus the old notion that increasing the order of 
spatial discretization leads to greater instability is dispelled if the time-march order is 
increased commensurately. Also, this increase in stability 
range is surprisingly linear with increase in
time-marching order. Each order gains $\approx 0.17$ in $\nu_c$. Thus the stability
range doubles in going from the first to the fourth-order.

\section {Solving nonlinear equations}

Nonlinear equations are difficult to solve in general. However, in the
case of the general nonlinear advection equation,
\be
\partial_t u=f(u)\partial_x u,
\la{nlin}
\ee
a simple formal solution exists and can be used to derive time-marching
algorithms of any order. By Taylor's expansion, one has
\be
u(x,\dt)=u(x,0)+\dt\partial_tu +\frac12 \dt^2\partial_t^2u +\frac1{3!} \dt^3\partial_t^3u+\cdots,
\ee 
where all the time-derivatives are evaluated at $t=0$. These derivatives can be obtained by
multiply both sides of (\ref{nlin}) by $f^j(u)$,
\ba
f^j(u)\partial_tu&=&f^{j+1}(u)\partial_x u,\nn\\
\quad\Rightarrow\quad \partial_t u_j(u)&=&\partial_x u_{j+1}(u),
\la{seq}
\ea
where we have defined, for $j\ge 0$,
\be
u_j(u)=	\int f^j(u)du,
\la{nlfct}
\ee
with $u_0(u)\equiv u(x,t)$. These are the conserved densities, since
\be
\partial_t\int_a^b u_j(u) dx=\int_a^b \partial_x u_{j+1}(u) dx=0
\ee
for periodic or Dirichlet boundary conditions. It follows from (\ref{seq}) that
\ba
\partial_tu&=&\partial_x u_1\nn\\
\partial_t^2u&=&\partial_x(\partial_t u_1 )
=\partial^2_x u_2\nn\\
&& \cdots\nn\\
\partial_t^j u&=&\partial^{j-1}_x(\partial_t u_{j-1}) =\partial^j_x u_j
\ea 
and therefore the solution is simply
\be
u(x,\dt)=u(x,0)+\dt\partial_x u_1 +\frac12 \dt^2\partial_x^2 u_2
+\frac1{3!} \dt^3\partial_x^3u_3 +\cdots.
\la{noneq}
\ee
To see how this solution works, consider the
inviscid Burgers' equation with 
$$f(u)=-u.$$
In this case, 
\be
u_n(x,t)=(-1)^n\frac{u^{n+1}(x,t)}{n+1}. 
\la{undef}
\ee
For the initial profile
\be
u(x,0)=u_0(x)\equiv
\begin{cases}
1   &{\rm if\ } x<0 \\
1-x &{\rm if\ } 0\le x \le 1\\
0   &{\rm if\  } x>0,
\end{cases}
\la{initp}
\ee
\be
\partial^n_x u_n=n!(1-x),
\ee
and the solution (\ref{noneq}) gives, for $1\ge u(x,\dt)\ge 0$,
\ba
u(x,\dt)&=&(1+\dt+\dt^2+\dt^3+\cdots)(1-x),\\
&=&\frac{(1-x)}{(1-\dt)}.
\ea
Thus the top edge of the wave at $u(x,\dt)=1$ moves with
unit speed, $x=\dt$, and the shock-front forms at $\dt=1$.
This formal solution is incapable of describing the motion of
the shock-front beyond $\dt=1$, but remarkably, as will be shown
below, finite-difference schemes base on it can. 

\begin{figure}[hbt]
\includegraphics[width=0.80\linewidth]{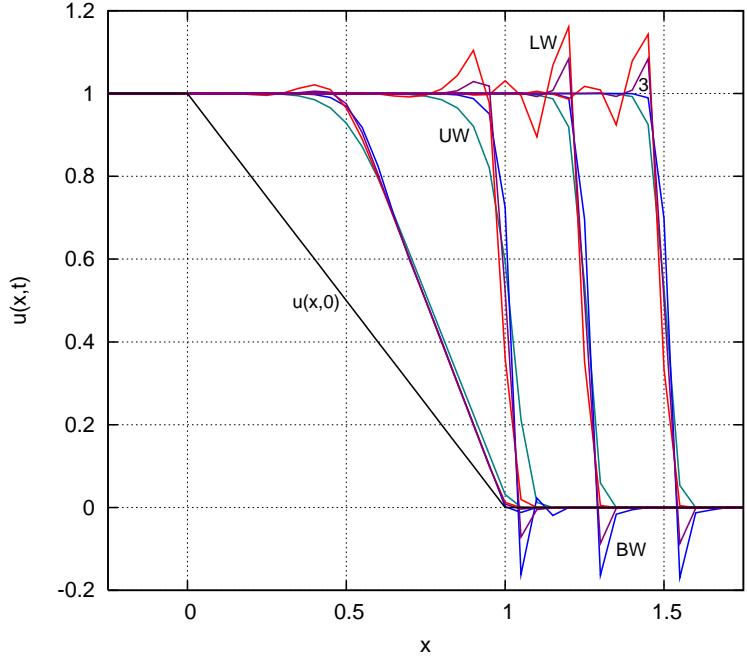}
\caption[]{\label{algad} 
The propagation of the inviscid Burgers' equation with initial profile (\ref{initp})
inside a [-5,5] box with $\dx=0.05$ and $\dt=0.025$. UW, LW, BW, and 3 denote the first order 
upwind, the second-order Lax-Wendroff, the second-order Beam-Warming, and the third-order 
algorithm described in the text.
The evolving profiles are given at $t=0,0.5,1.0,1.5$ and 2.0.
}
\end{figure}
 
The solution (\ref{noneq}) suggests that one should generalize the finite-difference scheme to
\be
u(x,\dt)=\sum_{i=1}^N c_{0i}\e^{k_i\dx\partial x} u(x,0)
	+\sum_{i=1}^N c_{1i}\e^{k_i\dx\partial x}u_1
	+\sum_{i=1}^N c_{2i}\e^{k_i\dx\partial x}u_2+\cdots.
\ee
Comparing this to (\ref{noneq}), one sees that
an $n$th-order time-marching algorithm now requires, in addition to $N=n+1$ grid points,
also nonlinear functions of $u(x,0)$ up to  $u_N(x,0)$.
For each $u_j(x,0)$, the set of coefficients $\{c_{ji}\}$ 
must have vanishing sums over all powers of $\{k_i\}$ up to $n$
except the following:
\be
\sum_{i=1}^N c_{0i}=1,\quad
\sum_{i=1}^N c_{1i}k_i=\nu, \quad
\sum_{i=1}^N c_{2i}k_i^2=\nu^2,\quad{\rm etc..}
\ee
where here $\nu=\dt/\dx$. Recalling (\ref{leq}), the solutions
are just
\be
c_{0i}=L_i(0),\quad
c_{1i}=\nu L^{(1)}_i(0), \quad
c_{2i}=\frac{\nu^2}{2!} L^{(2)}_i(0),\quad {\rm etc..}
\ee
and therefore the finite-difference scheme for solving (\ref{nlin}) is
\ba
u(x,\dt)
&=&\sum_{i=1}^{N}L_i(0)u(x+k_i\dx,0)+\nu \sum_{i=1}^{N}L_i^{(1)}(0)u_1(x+k_i\dx,0)\nn\\
&&+\frac{\nu^2}{2!}
\sum_{i=1}^{N}L_i^{(2)}(0)u_2(x+k_i\dx,0)
 +\frac{\nu^3}{3!} \sum_{i=1}^{N}L_i^{(3)}(0)u_3(x+k_i\dx,0) +\cdots\nn\\
\la{noeq}
\ea
If one were to replace all $u_n(x,0)\rightarrow u(x,0)$, then the above is just 
the linear advection scheme (\ref{anyad}) with $a_1=1$. Conversely,
{\it any linear advection scheme can now be used to solve the nonlinear 
advection equation} by replacing the $u(x,0)$ terms multiplying $\nu_1^n$ by $u_{n}(x,0)$. 
For example, the third-order advection scheme (\ref{ad3}) for solving the
inviscid Burgers' equation can now be applied here as
\ba
u^{k+1}_j
=u^k_j&+&\nu(\frac1{6}(u_1)^k_{j-2}-(u_1)^k_{j-1}+\frac12(u_1)^k_{j}+\frac13 (u_1)^k_{j+1})\nn\\
&+&\frac{\nu^2}{2!}((u_2)^k_{j-1}-2(u_2)^k_{j}+(u_2)^k_{j+1})\nn\\
&+&\frac{\nu^3}{3!}(-(u_3)^k_{j-2}+3(u_3)^k_{j-1}-3(u_3)^k_{j}+ (u_3)^k_{j+1})
\la{nlad3}
\ea
with $(u_n)^k_j=(-1)^n(u^k_j)^{n+1}/(n+1)$.
Thus arbitrary order schemes for solving the nonlinear advection equation (\ref{nlin}) can be obtained from
(\ref{anyad}).

To see how these schemes work, we compare their results when propagating the 
initial profile (\ref{initp}) from $t=0$ to $t=2$. Before the formation of the
shock front at $t=1$, the top edge of the wave is traveling at unit speed and reaches
$x=0.5$ and $x=1.0$ at $t=0.5$ and $t=1.0$ respectively. After the shock has formed,
the shock front travels at half the initial speed and reaches $x=1.25$ at $t=1.5$ and $x=1.5$
at $t=2.0$. The upwind (UW) scheme is overly diffusive, the Lax-Wendroff (LW) and the Beam-Warming (BW)
schemes have unwanted oscillations {\it trailing} and {\it ahead} of the shock front
respectively. Algorithm 3 has {\it reduced} oscillations both before and after the shock front.
While algorithms of any order for solving the linear advection equation is easily
generated, it remains difficult to produce arbitrary order algorithms for solving the
nonlinear advection equation, because one must disentangle each power of $\nu$ in
$c_i$ by hands.     

This example shows that for solving nonlinear equations, one must also discretize suitable
nonlinear functions of the propagating wave. For the nonlinear advection equation, the set of needed
nonlinear functions are given by (\ref{nlfct}). Unfortunately, the solution to the nonlinear 
{\it diffusion} equation is not of the form of (\ref{noneq}) and further study 
is needed to derive finite-difference schemes that can match its solution.
 
\section {Conclusions and discussions}
																	  
    In this work, we have shown that by matching the operator form of
the finite-difference scheme to the formal operator solution, one can 
systematically derive explicit finite-difference schemes for solving
any linear partial differential equation with constant coefficients.
This theory provided a unified description of all explicit finite-difference
schemes through the use of Lagrange polynomials. In a way, this work showed
that, not only are Lagrange polynomials important for doing interpolations, 
they are also cornerstones for deriving finite-difference schemes.

   Because one has a unified description of all finite-difference schemes,
there is no need to analyze each finite-difference scheme one by one. Theorem \ref{theerr},
for example, gives the local error for all algorithms at once. Also, 
the stability of first-order algorithms for solving (\ref{meq}) 
can be determined for all $m$ simultaneously. 
It would be of great interest if all second-order time-marching algorithms for solving 
the $m$-order linear equations can also be characterized the same way. 
The method used here for solving the operator equality (\ref{opeq}) is just 
Taylor's expansion, alternative methods of solving the equality without Taylor's expansion
would yield entirely new classes of finite-difference schemes.
												   
   Finally, this work focuses attention on obtaing the formal operator solution to 
the partial differential equation. 
To the extent that the formal solution embodies all the conservative properties of 
the equation, a sufficiently high-order approximation
to the formal solution should yield increasing better conservative schemes.
The method is surprisingly effective in deriving arbitrary order schemes for solving 
the general nonlinear advection equation (\ref{nlin}). One is therefore encouraged
to gain a deeper understanding of formal solutions so that better numerical schemes
can be derived for solving nonlinear equations.

\section*{Acknowledgment}

 This work is supported in part by a grant from the Qatar National Research 
 Fund \# NPRP 5-674-1-114.

 \Appendix
\section {Lagrange interpolation polynomials}
Consider the Lagrange interpolation at $N$ points $\{k_1, k_2, \dots, k_N\}$
with values $\{f_1, f_2, \dots, f_N\}$. The interpolating $N-1$ degree polynomial is given
by
\be
f(x)=\sum_{i=1}^N f_i L_i(x),
\la{lpoly}
\ee
where $L_i(x)$ are the Lagrange polynomials defined by
\be
L_i(x)=\prod_{j=1 (\ne i)}^N \frac{(x-k_j)}{(k_i-k_j)}.
\ee
Since by construction
\be
L_i(k_j)=\delta_{ij}
\ee
one has the desired interpolation,
\be
f(k_j)=\sum_{i=1}^n f_i L_i(k_j)=\sum_{i=1}^n f_i \delta_{ij}=f_j.
\ee
Now let $f_i=k_i^m$ for $0\le m\le N-1$, 
then the interpolating polynomial
\be
f(x)=\sum_{i=1}^n k_i^m L_i(x)
\ee
and the function
\be
g(x)=x^m
\ee
both interpolate the same set of points and therefore must agree. Hence
\be
\sum_{i=1}^N k^m_iL_i(x)=x^m,
\ee
for $0\le m\le N-1$.
This is then (\ref{solut2}).



\end{document}